\magnification 1200
\hsize 13.2cm
\vsize 20cm
\parskip 3pt plus 1pt
\parindent 5mm

\def\\{\hfil\break}


\font\seventeenbf=cmbx10 at 17.28pt

\font\twelvebf=cmbx10 at 12pt
\font\eightbf=cmbx8
\font\sixbf=cmbx6

\font\eighti=cmmi8
\font\sixi=cmmi6

\font\eightrm=cmr8
\font\sixrm=cmr6

\font\eightsy=cmsy8
\font\sixsy=cmsy6

\font\eightit=cmti8
\font\eighttt=cmtt8
\font\eightsl=cmsl8

\font\seventeenbsy=cmbsy10 at 17.28pt

\font\twelvebsy=cmbsy10 at 12pt
\font\tenbsy=cmbsy10
\font\eightbsy=cmbsy8
\font\sevenbsy=cmbsy7
\font\sixbsy=cmbsy6
\font\fivebsy=cmbsy5

\font\tenmsa=msam10

\font\sevenmsa=msam7
\font\fivemsa=msam5
\newfam\msafam
  \textfont\msafam=\tenmsa
  \scriptfont\msafam=\sevenmsa
  \scriptscriptfont\msafam=\fivemsa

\font\tenmsb=msbm10
\font\eightmsb=msbm8
\font\sevenmsb=msbm7
\font\fivemsb=msbm5
\newfam\msbfam
  \textfont\msbfam=\tenmsb
  \scriptfont\msbfam=\sevenmsb
  \scriptscriptfont\msbfam=\fivemsb
\def\Bbb{\fam\msbfam\tenmsb}

\font\tenCal=eusm10
\font\sevenCal=eusm7
\font\fiveCal=eusm5
\newfam\Calfam
  \textfont\Calfam=\tenCal
  \scriptfont\Calfam=\sevenCal
  \scriptscriptfont\Calfam=\fiveCal
\def\Cal{\fam\Calfam\tenCal}

\font\teneuf=eusm10
\font\teneuf=eufm10
\font\seveneuf=eufm7
\font\fiveeuf=eufm5
\newfam\euffam
  \textfont\euffam=\teneuf
  \scriptfont\euffam=\seveneuf
  \scriptscriptfont\euffam=\fiveeuf

\font\seventeenbfit=cmmib10 at 17.28pt

\font\twelvebfit=cmmib10 at 12pt
\font\tenbfit=cmmib10
\font\eightbfit=cmmib8
\font\sevenbfit=cmmib7
\font\sixbfit=cmmib6
\font\fivebfit=cmmib5
\newfam\bfitfam
  \textfont\bfitfam=\tenbfit
  \scriptfont\bfitfam=\sevenbfit
  \scriptscriptfont\bfitfam=\fivebfit


\catcode`\@=11
\def\eightpoint{%
  \textfont0=\eightrm \scriptfont0=\sixrm \scriptscriptfont0=\fiverm
  \def\rm{\fam\z@\eightrm}%
  \textfont1=\eighti \scriptfont1=\sixi \scriptscriptfont1=\fivei
  \def\oldstyle{\fam\@ne\eighti}%
  \textfont2=\eightsy \scriptfont2=\sixsy \scriptscriptfont2=\fivesy
  \textfont\itfam=\eightit
  \def\it{\fam\itfam\eightit}%
  \textfont\slfam=\eightsl
  \def\sl{\fam\slfam\eightsl}%
  \textfont\bffam=\eightbf \scriptfont\bffam=\sixbf
  \scriptscriptfont\bffam=\fivebf
  \def\bf{\fam\bffam\eightbf}%
  \textfont\ttfam=\eighttt
  \def\tt{\fam\ttfam\eighttt}%
  \textfont\msbfam=\eightmsb
  \def\Bbb{\fam\msbfam\eightmsb}%
  \abovedisplayskip=9pt plus 2pt minus 6pt
  \abovedisplayshortskip=0pt plus 2pt
  \belowdisplayskip=9pt plus 2pt minus 6pt
  \belowdisplayshortskip=5pt plus 2pt minus 3pt
  \smallskipamount=2pt plus 1pt minus 1pt
  \medskipamount=4pt plus 2pt minus 1pt
  \bigskipamount=9pt plus 3pt minus 3pt
  \normalbaselineskip=9pt
  \setbox\strutbox=\hbox{\vrule height7pt depth2pt width0pt}%
  \let\bigf@ntpc=\eightrm \let\smallf@ntpc=\sixrm
  \normalbaselines\rm}
\catcode`\@=12

\def\eightpointbf{%
 \textfont0=\eightbf   \scriptfont0=\sixbf   \scriptscriptfont0=\fivebf
 \textfont1=\eightbfit \scriptfont1=\sixbfit \scriptscriptfont1=\fivebfit
 \textfont2=\eightbsy  \scriptfont2=\sixbsy  \scriptscriptfont2=\fivebsy
 \eightbf
 \baselineskip=10pt}

\def\tenpointbf{%
 \textfont0=\tenbf   \scriptfont0=\sevenbf   \scriptscriptfont0=\fivebf
 \textfont1=\tenbfit \scriptfont1=\sevenbfit \scriptscriptfont1=\fivebfit
 \textfont2=\tenbsy  \scriptfont2=\sevenbsy  \scriptscriptfont2=\fivebsy
 \tenbf}

\def\twelvepointbf{%
 \textfont0=\twelvebf   \scriptfont0=\eightbf   \scriptscriptfont0=\sixbf
 \textfont1=\twelvebfit \scriptfont1=\eightbfit \scriptscriptfont1=\sixbfit
 \textfont2=\twelvebsy  \scriptfont2=\eightbsy  \scriptscriptfont2=\sixbsy
 \twelvebf
 \baselineskip=14.4pt}

\def\seventeenpointbf{%
 \textfont0=\seventeenbf  \scriptfont0=\twelvebf  \scriptscriptfont0=\eightbf
 \textfont1=\seventeenbfit\scriptfont1=\twelvebfit\scriptscriptfont1=\eightbfit
 \textfont2=\seventeenbsy \scriptfont2=\twelvebsy \scriptscriptfont2=\eightbsy
 \seventeenbf
 \baselineskip=20.736pt}


\newdimen\srdim \srdim=\hsize
\newdimen\irdim \irdim=\hsize
\def\NOSECTREF#1{\noindent\hbox to \srdim{\null\dotfill ???(#1)}}
\def\SECTREF#1{\noindent\hbox to \srdim{\csname REF\romannumeral#1\endcsname}}
\def\INDREF#1{\noindent\hbox to \irdim{\csname IND\romannumeral#1\endcsname}}
\newlinechar=`\^^J
\def\openauxfile{
  \immediate\openin1\jobname.aux
  \ifeof1
  \message{^^JCAUTION\string: you MUST run TeX a second time^^J}
  \let\sectref=\NOSECTREF \let\indref=\NOSECTREF
  \else
  \input \jobname.aux
  \message{^^JCAUTION\string: if the file has just been modified you may
    have to run TeX twice^^J}
  \let\sectref=\SECTREF \let\indref=\INDREF
  \fi
  \message{to get correct page numbers displayed in Contents or Index
    Tables^^J}
  \immediate\openout1=\jobname.aux
  \let\END=\end \def\end{\immediate\closeout1\END}}

\newbox\titlebox   \setbox\titlebox\hbox{\hfil}
\newbox\sectionbox \setbox\sectionbox\hbox{\hfil}
\def\folio{\ifnum\pageno=1 \hfil \else \ifodd\pageno
           \hfil {\eightpoint\copy\sectionbox\kern8mm\number\pageno}\else
           {\eightpoint\number\pageno\kern8mm\copy\titlebox}\hfil \fi\fi}
\footline={\hfil}
\headline={\folio}

\def\titlerunning#1{\setbox\titlebox\hbox{\eightpoint #1}}
\def\title#1{\noindent\hfil$\smash{\hbox{\seventeenpointbf #1}}$\hfil
             \titlerunning{#1}\medskip}

\newcount\numbersection \numbersection=-1
\def\sectionrunning#1{\setbox\sectionbox\hbox{\eightpoint #1}
  \immediate\write1{\string\def \string\REF
      \romannumeral\numbersection \string{%
      \noexpand#1 \string\dotfill \space \number\pageno \string}}}
\def\section#1{%
  \par\vskip0.666cm\penalty -100
  \vbox{\baselineskip=14.4pt\noindent{{\twelvepointbf #1}}}
  \vskip2pt
  \penalty 500
  \advance\numbersection by 1
  \sectionrunning{#1}}

\def\subsection#1{%
  \par\vskip0.5cm\penalty -100
  \vbox{\noindent{{\tenpointbf #1}}}
  \vskip1pt
  \penalty 500}

\newcount\numberindex \numberindex=0
\def\index#1#2{%
  \advance\numberindex by 1
  \immediate\write1{\string\def \string\IND #1%
     \romannumeral\numberindex \string{%
     \noexpand#2 \string\dotfill \space \string\S \number\numbersection,
     p.\string\ \space\number\pageno \string}}}

\newdimen\itemindent \itemindent=\parindent

\def\item#1{\par\noindent\hangindent\itemindent%
            \rlap{#1}\kern\itemindent\ignorespaces}
\def\itemitem#1{\par\noindent\hangindent2\itemindent%
            \kern\itemindent\rlap{#1}\kern\itemindent\ignorespaces}
\def\itemitemitem#1{\par\noindent\hangindent3\itemindent%
            \kern2\itemindent\rlap{#1}\kern\itemindent\ignorespaces}

\long\def\claim#1|#2\endclaim{\par\vskip 5pt\noindent
{\tenpointbf #1.}\ {\it #2}\par\vskip 5pt}

\def\today{\ifcase\month\or
January\or February\or March\or April\or May\or June\or July\or August\or
September\or October\or November\or December\fi \space\number\day,
\number\year}

\catcode`\@=11
\newcount\@tempcnta \newcount\@tempcntb
\def\timeofday{{%
\@tempcnta=\time \divide\@tempcnta by 60 \@tempcntb=\@tempcnta
\multiply\@tempcntb by -60 \advance\@tempcntb by \time
\ifnum\@tempcntb > 9 \number\@tempcnta:\number\@tempcntb
  \else\number\@tempcnta:0\number\@tempcntb\fi}}
\catcode`\@=12

\def\bibitem#1&#2&#3&#4&%
{\hangindent=0.8cm\hangafter=1
\noindent\rlap{\hbox{\eightpointbf #1}}\kern0.8cm{\rm #2}{\it #3}{\rm #4.}}


\def\bC{{\Bbb C}}

\def\bR{{\Bbb R}}


\def\cC{{\Cal C}}
\def\cD{{\Cal D}}


\def\square{{\hfill \hbox{
\vrule height 1.453ex  width 0.093ex  depth 0ex
\vrule height 1.5ex  width 1.3ex  depth -1.407ex\kern-0.1ex
\vrule height 1.453ex  width 0.093ex  depth 0ex\kern-1.35ex
\vrule height 0.093ex  width 1.3ex  depth 0ex}}}
\def\qed{\kern10pt$\square$}
\def\hexnbr#1{\ifnum#1<10 \number#1\else
 \ifnum#1=10 A\else\ifnum#1=11 B\else\ifnum#1=12 C\else
 \ifnum#1=13 D\else\ifnum#1=14 E\else\ifnum#1=15 F\fi\fi\fi\fi\fi\fi\fi}
\def\msatype{\hexnbr\msafam}
\def\msbtype{\hexnbr\msbfam}
\mathchardef\restriction="3\msatype16   
\mathchardef\boxsquare="3\msatype03
\mathchardef\preccurlyeq="3\msatype34
\mathchardef\compact="3\msatype62
\mathchardef\smallsetminus="2\msbtype72   
\mathchardef\subsetneq="3\msbtype28
\mathchardef\supsetneq="3\msbtype29
\mathchardef\leqslant="3\msatype36   
\mathchardef\geqslant="3\msatype3E   
\mathchardef\stimes="2\msatype02
\mathchardef\ltimes="2\msbtype6E
\mathchardef\rtimes="2\msbtype6F

\def\dbar{\overline\partial}
\def\ddbar{\partial\overline\partial}

\let\ol=\overline

\let\wt=\widetilde
\let\wh=\widehat
\let\text=\hbox
\def\buildo#1^#2{\mathop{#1}\limits^{#2}}
\def\buildu#1_#2{\mathop{#1}\limits_{#2}}
\def\ort{\mathop{\hbox{\kern1pt\vrule width4.0pt height0.4pt depth0pt
    \vrule width0.4pt height6.0pt depth0pt\kern3.5pt}}}

\def\vlra{\mathrel{\smash-}\joinrel\mathrel{\smash-}\joinrel%
    \kern-2pt\longrightarrow}
\def\srelbar{\vrule width0.6ex height0.65ex depth-0.55ex}
\def\merto{\mathrel{\srelbar\kern1.3pt\srelbar\kern1.3pt\srelbar
    \kern1.3pt\srelbar\kern-1ex\raise0.28ex\hbox{${\scriptscriptstyle>}$}}}

\def\Re{\mathop{\rm Re}\nolimits}

\def\Hom{\mathop{\rm Hom}\nolimits}

\def\deg{\mathop{\rm deg}\nolimits}


\long\def\InsertFig#1 #2 #3 #4\EndFig{\par
\hbox{\hskip #1mm$\vbox to#2mm{\vfil\special{"
(/home/demailly/psinputs/grlib.ps) run
#3}}#4$}}
\long\def\LabelTeX#1 #2 #3\ELTX{\rlap{\kern#1mm\raise#2mm\hbox{#3}}}


\itemindent = 4mm

\title{An Ohsawa-Takegoshi theorem on}
\title{compact K\"{a}hler manifolds}


\vskip10pt

\centerline{Li YI}

\centerline{Institut Elie Cartan, Nancy}

\vskip20pt

\noindent{\bf Abstract. \it {In this article we prove a
theorem of Ohsawa-Takegoshi type on compact K\"{a}hler manifolds.
Our arguments follow the ``standard" approach for this kind of extension results; however,  
there are many complications arising from the regularization process
of quasi-psh functions on compact K\"ahler manifolds, and unfortunately
we only obtain a particular case of the expected result. We remark that the additional hypothesis we are forced to make are natural, since they are 
verified in many situations; we hope to remove them in a near future.}}

%
\vskip 15pt 
\section{\S0 Introduction}
\vskip 5pt 

\noindent We begin by fixing a few notations, and by presenting the 
general context in which we will work.

Let $X$ be a compact K\"{a}hler manifold and let $Z\subset X$ be the
zero set of a holomorphic section $s\in H^{0}(X,E)$ of a hermitian
line bundle $(E,h_{E})$. The hypersurface $Z$ is assumed to be
non-singular. Let $L$ be a line bundle,
endowed with a possibly singular metric $h_{L}$ such that:

{\itemindent 0.5cm

\item{(1)} $\displaystyle \Theta_{h_{L}}(L)\geq 0$ as a current on $X$;
\smallskip
\item{(2)} $\displaystyle \Theta _{h_{L}}(L) \geq
{1\over \alpha}\Theta(E)$ for some
$\alpha\geq 1, \quad \alpha\in{\cal C}^\infty(X)$;
\smallskip
\item{(3)} $\displaystyle|s|^{2}_{h_{E}}\leq \exp(-\alpha)$ on $X$, and
the restriction of the metric $h_{L}$ to $Z$ is well defined.}
\medskip

\noindent We recall next the notion of multiplier ideal sheaf 
associated to the metric $h_L$; it is a way of measuring the ``large"
algebraic singularities this metric may have, and it plays an important role in
both algebraic and analytic geometry. Hence we introduce 
${\cal I}_{\varphi_{L}}$, defined locally at each point $x\in X$ as follows

$${\cal
I}({\varphi_{L}})_x:= \{f\in {\cal O}_{X, x}: |f|^2e^{-\varphi_L}\in L^1(X, x)\}.
$$

Actually, a version of this ideal (introduced and investigated in [5]) will be relevant for us in what follows; we define
$${\cal
I}_+({\varphi_{L}}):= \lim_{\varepsilon\to 0}
{\cal I}\big((1+ \varepsilon){\varphi_{L}}\big).$$
We notice that we have 
$${\cal
I}_+({\varphi_{L}})= {\cal I}\big((1+ \varepsilon_0){\varphi_{L}}\big)$$
for some positive real number $\varepsilon_0$.
Conjecturally (cf. [5]) we have ${\cal
I}({\varphi_{L}})= {\cal
I}_+({\varphi_{L}})$, but so far this equality is only established up to the 
dimension 2 (and it seems to be a difficult problem).

\medskip

\noindent We will formulate next a requirement concerning the singularities of
the metric $\varphi_L$ locally near the  set $Z$; we note that this condition equally appears in [15].
\smallskip

\noindent {\bf Hypothesis ${\cal C}$.} {\sl Let $x_0\in Z$ be an arbitrary point,
and let $z= (z_1,..., z_n)$ be a coordinate system of $X$ defined on an open set $U$, such that $\displaystyle Z\cap U= (z_1= 0)$. Let $\displaystyle f\in 
{\cal I}({\varphi_{L}|_Z})_{x_0}$ be an element of the multiplier sheaf associated to $\varphi_{L}|_Z$ defined on $Z\cap U$, and let $\widetilde f$ be an extension of $f$ to $U$. Then the function
$$t\to \int_{U^\prime}|\widetilde f(t, z^\prime)|^2e^{-\varphi_L(t, z^\prime)}d\lambda(z^\prime)$$
is continuous at $t= 0$, where $z^\prime:= (z^2,..., z^n)$.
}

\medskip

\noindent We state our main theorem as follows.

\claim 0.1 Theorem|Let $X$ be a compact K\"ahler manifold, 
and let $(L, h_L)$ be a hermitian line bundle, such that the properties {\rm (1)-(3)} above
together with the hypothesis $\cal C$ are satisfied. 
Given a section 
$$u\in H^{0}\big(Z,(K_{Z}+L_{|Z})\otimes 
{\cal I}_+(h_{L}|_{Z})\big)$$ such that
there exists a holomorphic section $U$ of $K_X+ Z+ L$ such that
$U|_Z=u\wedge ds$ and such that 
$$\int_{X}{|U|^2 e^{-\varphi_{L}-\varphi_E}\over |s|^{2}\log^2(|s|^2)}\leq
 C_0\int_{Z}|u|^{2} e^{-\varphi_{L}}$$
where $C_0$ is a purely numerical constant.
\endclaim
\vskip 30pt

\noindent As we have already mentioned in the abstract, our proof follows the 
general outline of the classical works on the subject (starting with [12], [13], [14], [4], [7], [8], [10], [17], [18], [19]).
The additional -and rather severe- difficulty comes from the procedure of 
regularization of quasi-psh functions on compact complex manifolds. To put this in a proper perspective, we recall here in a few words the approach in the projective case.

Let $X$ be a {\sl projective} manifold; then there exists a hyperplane section 
$H$ which do not contain $Z$, and such that 
$$X\setminus H= \bigcup_k \Omega_k$$
where $(\Omega_k)$ is an increasing sequence of Stein sub-domains of $X$.
The restriction of $h_L$ to each $\Omega_k$ can be written as (decreasing) limit of $\cC^\infty$ metrics $h_{L, \rho}$, in such a way that the curvature/normalization hypothesis in the main theorem above are preserved. This is precisely the reason for which the Stein property of $\Omega_k$ is used: we can assume that $\displaystyle h_L|_{\Omega_k}$ is non-singular, provided that the estimates we obtain for the norm of the extension are independent of $k$.
In conclusion, we use the well-known version of Ohsawa-Takegoshi theorem on
each $\Omega_k$, and then let $k\to\infty$. We remark here that there are 3 parameters involved in the proof: $k$, $\varepsilon$ (used to concentrate the mass on $Z\cap\Omega_k$, see [1], [17]) and $\rho$: they are completely independent,
i.e. for each $k$, we are allow to first let $\varepsilon\to 0$, and then $\delta\to 0$.

By contrast, if $X$ is only assumed to be K\"ahler, then we do not have at our 
disposal the sequence $(\Omega_k)$. However, we can still regularize
the metric $h_L$ and obtain a sequence $h_{L, \rho}$ by using the main result in [3] (in this way we can compensate the absence of the ``large" coordinate systems used in the projective case).
The important difference is that in the actual context, the curvature hypothesis
(1) and (2) above are only verified up to a negative factor $-\delta_\rho\omega$,
where $\delta_\rho\to 0$ as $\rho\to 0$. The occurrence of this small negative factor has a huge
consequence: we cannot allow $\varepsilon$ and $\rho$ to be independent
anymore, mainly because of the estimates we obtain as a consequence of the twisted
Bochner identity. Basically, this is the reason why we are forced to make
the additional hypothesis  in the main theorem: they are needed in order to insure that a few limit processes involved in our proof are justified
in the context of the dependence of the parameters $\rho, \varepsilon$.

We have divided our arguments in several steps. In the first part of the proof, we will use an approximation theorem for quasi-psh functions,
due to J.-P. Demailly (cf. [3]).
Then we construct a smooth extension of the given section $u$, and we analyze its properties (norm of its $\dbar$, etc). Finally, we convert this extension to a holomorphic one, by solving a $\dbar$ equation, using a twisted Bochner inequality.\\

{\it Acknowledgements.}  I would like to thank sincerely my supervisor Mihai Paun for introducing me this interesting problem, for many useful discussions on the difficulties I encountered while writing this article and, finally, for helping me polishing my mathematical writing.  

\section{\S A Approximation of quasi-psh functions}

\vskip 5pt

\noindent 
The regularization result we need is the following (cf. [3], page 18).

\claim Theorem ([3])|Let $T:= \alpha+ \sqrt{-1}\ddbar \varphi$ be a closed current on a compact complex manifold $X$; here $\alpha$ is a closed 
(1,1)-form. We assume that $\varphi$ is quasi-psh, and let $\gamma$ be a 
continuous form of $(1, 1)$-type, such that $T\geq \gamma$.
Then for each $\rho> 0$, there exists a function $\varphi_\rho\in L^1(X)$ such that:

{\itemindent 0.6cm

\item {(i)} The function $\varphi_\rho$ is smooth on $X\setminus E_\rho (T)$, where $E_\rho(T)$ is the $\rho$-upper level set
of Lelong numbers of $T$.

\item {(ii)} We have $T_\rho:= \alpha+ \sqrt{-1}\ddbar \varphi_\rho\geq\gamma-\delta_\rho\omega$, where $\delta_\rho\to 0$ as $\rho\to 0$.

\item{(iii)} We have $\varphi\leq \varphi_\rho$, for each $\rho > 0$.
}

\endclaim

\medskip

\noindent One of the important aspects of the previous result is that the construction of $\varphi_\rho$ is independent of $\gamma$, in the following sense. Given a quasi-psh function $\varphi$ on $X$, the family of
functions $(\varphi_\rho)_{\rho> 0}$ is obtained by using a regularization kernel, followed by a Legendre transform. As we can see (cf. [3]),
these operations are independent of the form $\gamma$; in conclusion, given two forms $\gamma_1$ and $\gamma_2$ such that 
$$T\geq \gamma_1, \quad T\geq \gamma_2$$ 
we have 

$$T_\rho\geq \gamma_j- \delta_\rho\omega$$
for each $j= 1, 2$.

\noindent Concretely, the construction is done as follows: we first construct  a modified exponential map
$$T_X\to X \quad (x,\zeta)\to \text{exph}_x(\zeta) \quad \zeta\in T_{X,x}$$
from the tangent bundle of $X$ to $X$ with the property that for every $x\in X$, the map $\zeta\to \text{exph}_x(\zeta)$ has a holomorphic Taylor expansion at $\zeta=0.$ Then we set $\varphi^\prime_\rho(z)$ to be the convolution of $\varphi(\text{exph}_x(\zeta))$ with a cut-off function $\chi:  \Bbb R\to \Bbb R$ of class $\cal C^\infty:$ 
$$\varphi^\prime_\rho(z)={1\over \rho^{2n}}\int_{\zeta\in T_{X,z}}\varphi(\text{exph}_z(\zeta))\chi\left({|\zeta|^2 \over \rho^2}\right)d\lambda(\zeta), \quad \rho>0.$$
Here the cut-off function $\chi$ satisfies that
$$\chi(t)>0\quad\text{for}\quad t<1,\quad\chi(t)=0\quad\text{for}\quad t\geq 1,\quad\int_{v\in{\Bbb C}^n}\chi(|v|^2)d\lambda(v)=1$$
and $d\lambda$ denotes the Lebesgue measure on the hermitian space $(T_{X,z},\omega(z))$ resp. on ${\Bbb C}^n.$ Set then $\Phi(z,w)=\varphi^\prime_\rho(z)$ for $w\in \Bbb C,$ $|w|=\rho$ with
$$\Phi(z,w)=\int_{\zeta\in T_{X,z}}\varphi(\text{exph}_z(w\zeta))\chi(|\zeta|^2)d\lambda(\zeta).$$
If we change the variable $y=\text{exph}_z(w\zeta),$ we can see that $w\zeta$ is a smooth function of $y,z$ in a neighborhood of the diagonal in $X\times X.$ Hence $\Phi$ is smooth over $X\times \{0<|w|<\rho_0\}$ for some $\rho_0>0.$ The family of $\varphi_\rho$ are then defined by the Legendre transform 
$$\varphi_\rho=\inf\limits_{|w|<1}\left(\wt\Phi(x,\rho w)+{\rho\over{1-|w|^2}}-\rho\log|w|\right),$$
where $\wt\Phi(x,w)=\Phi(x,w)+|w|.$ So its construction is independent of $\gamma.$ We refer to [3] for more details.

\medskip

\noindent We apply these considerations to the current 
$$T:= \Theta_{h_L}(L)$$
and we formulate the conclusion as follows.

\claim Lemma|For each $\rho> 0$ there exists a metric 
$h_{L, \rho}= e^{-\varphi_{L, \rho}}$ on $L$ such that we have:
$$\Theta_{h_{L, \rho}}(L)\geq -\delta_\rho\omega, \quad
\Theta_{h_{L, \rho}}(L)\geq {1\over \alpha}\Theta(E)-\delta_\rho\omega.$$
Moreover, there exists an analytic set $Y_\rho\subset X$ such that 
$\varphi_{L, \rho}$ is smooth on $X\setminus Y_\rho$, for each 
$\rho> 0$.
\endclaim

\noindent Indeed, this is a consequence of the regularization theorem quoted above, together with the hypothesis (1) and (2) of the main theorem.


\section{\S B. Construction of a ${\cal C}^\infty$ extension and its properties}

\medskip 

\noindent We consider here coordinates sets isomorphic to polydisks
$U_\alpha\subset X$, together with local coordinates $z_\alpha= (z^1_\alpha,..., z^n_\alpha)$ such that $Z\cap U_\alpha= (z^1_\alpha= 0)$. Let 
$f_\alpha$ be the holomorphic function defined on $Z\cap U_\alpha$, which corresponds to the restriction of the global section $u$ to the previous set.
By hypothesis, we have $\displaystyle u\in {\cal I}_+(h_{L}|_Z)$ which in local terms
translates as
$$\int_{Z\cap U_\alpha}|f_\alpha|^2e^{-(1+ \varepsilon_0)\varphi_L}
d\lambda(z_\alpha)< \infty;$$
here $\varepsilon_0$ is a small, positive real number.

By the version of Ohsawa-Takegoshi theorem in [3], for each index $\alpha$ there exists a holomorphic function 
$\wt f_\alpha$ defined on $U_\alpha$, such that 
$\wt f_\alpha|_{Z\cap U_\alpha}= f_\alpha$, and such that 
$$\int_{U_\alpha}
{|\wt f_\alpha|^2e^{-(1+ \varepsilon_0)\varphi_L}\over |s|^2\log^2{1\over |s|^2}}d\lambda(z_\alpha)\leq 
C_0\int_{Z\cap U_\alpha}|f_\alpha|^2e^{-(1+ \varepsilon_0)\varphi_L}
d\lambda(z_\alpha)\leqno(1)$$ 
where $C_0$ above is a numerical constant.

Let $(\theta_\alpha)$ be a partition of unit corresponding to the covering $(U_\alpha)$; we define the section $U_\infty$ of $K_X+ Z+ L$ as follows
$$U_\infty:= \sum_\alpha\theta_\alpha \wt f_\alpha $$
where $\wt f_\alpha$ are seen as local holomorphic sections of the bundle
above, which are globalized via the partition of unit.

\noindent We will use $U_\infty$ in order to construct the holomorphic extension $U$ required by the main theorem, but prior to this, 
we analyze next the $L^2$ properties of $\dbar U_\infty$. To this end,
we remark that for any index $\beta$ we have
$$\eqalign{
\dbar U_{\infty}|_{U_\beta}= & \sum_\alpha \wt f_\alpha \dbar \theta_\alpha|_{U_\beta}= \cr
= & \sum_\alpha (\wt f_\alpha- \wt f_\beta) \dbar \theta_\alpha|_{U_\beta}. \cr
}
 \leqno (2)$$

Let us consider an index $\gamma$ such that $U_\gamma\cap U_\beta\neq\emptyset$. Then there exists a holomorphic function $g_{\gamma\beta}$
defined on $U_\gamma\cap U_\beta$ such that 
$$\wt f_\gamma- \wt f_\beta= s g_{\gamma\beta}\leqno(3)$$
since $\wt f_\gamma$ and $\wt f_\beta$ are both local extensions of $u$.
Moreover, the next relation holds
$$\int_{U_\gamma\cap U_\beta}{|g_{\gamma\beta}|^2e^{-(1+ \varepsilon_0)\varphi_L}\over \log^2{1\over |s|^2}}d\lambda< \infty$$
as a consequence of (3). 

We denote by $\varepsilon$ a positive number; during the proof of the main theorem, we will have to evaluate an expression of the following type 
$$I_\varepsilon:= \int_{|s|< \varepsilon}{|\dbar U_\infty|_\omega^2e^{-\varphi_L}\over |s|^2}dV_\omega;
\leqno(4)$$
as $\varepsilon\to 0$.  
We remark that up to a constant which is independent of $\varepsilon$, the quantity $I_\varepsilon$ is bounded from above by a sum of integrals of the following type
$$\int_{U_\alpha\cap U_\beta\cap |s|< \varepsilon} 
|g_{\alpha\beta}|^2e^{-\varphi_L};\leqno (5)$$
by H\"older inequality the quantity (4) is smaller than
$$\eqalign{
& \Big(\int_{U_\alpha\cap U_\beta}{|g_{\alpha\beta}|^2e^{-(1+ \varepsilon_0)\varphi_L}\over \log^2{1\over |s|^2}}d\lambda\Big)^{1\over 1+ \varepsilon_0}\Big(\int_{U_\alpha\cap U_\beta\cap |s|< \varepsilon}
|g_{\alpha\beta}|^2\log^{2\over \varepsilon_0}{1\over |s|^2}d\lambda\Big)^
{\varepsilon_0\over 1+\varepsilon_0}\leq \cr
\leq & C\big(\varepsilon^2\log^{2\over \varepsilon_0}{1\over \varepsilon^2}\big)^{\varepsilon_0\over 1+\varepsilon_0}.\cr
}$$

\noindent Hence we have
$$I_\varepsilon\leq C\big(\varepsilon^2\log^{2\over \varepsilon_0}{1\over \varepsilon^2}\big)^{\varepsilon_0\over 1+\varepsilon_0}\leqno(6)$$
where $C> 0$ is a constant independent of $\varepsilon$.
We remark that so far, we have only used the (extra) assumption 
$u\in {\cal I}_+(h_L)$, there are no restrictions on the singularities of 
$\varphi_L$ which are needed.
In the next paragraph, we follow the usual approach for the Ohsawa-Takegoshi type theorems by recalling the twisted Bochner formula in some detail, for the 
comfort of the reader.

\section{\S C. Twisted a-priori inequality}

\vskip 10pt Let $X$ be a complex manifold, endowed with a hermitian metric $\omega$. If 
$z^1,..., z^n$ are local coordinates on some open set $U\subset X$, then the metric $\omega$ viewed as a 
$(1,1)$-form can be expressed as
$$\omega= \sqrt {-1}\sum_{j, k}\omega_{j\ol k}dz^j\wedge dz^{\ol k}.$$
We are particularly interested in the tensor bundles $\displaystyle \Lambda ^{p, q}T^\star _X:= 
\Lambda ^{p}T^\star _X\otimes \Lambda ^{q}\ol T^\star _X$; we
can use $\omega$ to endow each of them with a hermitian structure (see e.g. [2], [11]).

More generally, let $E\to X$ be a holomorphic vector bundle of rank $r$ and let $h$ be a hermitian 
metric on $E$. Then we have a metric induced on the space of $(p,q)$-forms with values in $E$ with measurable coefficients and let us set 
$$\Vert u\Vert ^2:= \int_X\vert u\vert ^2dV_\omega;$$
we can equally define in a similar fashion a inner product
$$\langle\! \langle u, v\rangle\!\rangle:= \int_X\langle u, v\rangle dV_\omega.
$$
We obtain in this way a Hilbert space usually denoted by $L^2(X, \Lambda ^{p, q}T^\star _X\otimes E)$; the space of smooth $(p, q)$ forms with values in $E$ and compact support 
is a dense subspace of the above Hilbert space, and will be denoted with $\displaystyle \cD(X, \Lambda ^{p, q}T^\star _X\otimes E)$.

The holomorphic structure of $E$ allows us to extend the $\dbar$ to the space of 
sections of $E$; the reason is that the transition functions of $E$ belong to the kernel of 
$\dbar$. 

In a similar fashion, we can define the $\partial$ of a section (or more generally, of a
$(p,q)$ form) of a vector bundle whose transition functions are {\sl anti-holomorphic}.
This observation can be used to define an operator acting on the space of the $(p, q)$
forms with values in a {\it holomorphic } vector bundle which is the analog of the full
exterior derivative $d$.

Indeed, given a holomorphic hermitian bundle $(E, h)$, we have a metric identification 
$\rho: E\to \ol E^\star$, defined by $\rho(v)(\ol w):= \langle v, w\rangle _h$ and we can use this
identification in order to define the next differential operator:
$$\partial_h (u):= \rho^{-1}\partial \big(\rho(u)\big).$$
In local coordinates, the operator $D^\prime$ can be described as follows. If $(e_\alpha)_{\alpha=1,...,r}$ is a holomorphic local frame of $E_{|\Omega}$, let us define 
$$h_{\alpha \ol \beta}:= \langle e_\alpha, e_\beta\rangle.$$
The section $u$ can be written locally as follows $u= \sum_\alpha u^\alpha e_\alpha $ and then
$$\rho(u)= \sum u^\alpha h_{\alpha \ol \beta}e^{\star \ol \beta}$$
where the $(e^{\star \ol \beta})$ is the induced frame on $\ol E^\star$. Next we compute
$$\partial \rho(u)= \sum (\partial u^\alpha h_{\alpha \ol \beta}+ u^\alpha \partial h_{\alpha \ol \beta}) e^{\star \ol \beta}$$
and finally we obtain 
$$\eqalign{
\partial_h (u)= & \sum h^{\ol \beta\gamma }(\partial u^\alpha h_{\alpha \ol \beta}+ u^\alpha \partial h_{\alpha \ol \beta})e_\gamma= \cr
= & \sum (\partial u^\gamma+ u^\alpha \partial h_{\alpha \ol \beta}h^{\ol \beta\gamma })e_\gamma.
}$$
Of course, this is nothing but the $(1,0)$ part of the Chern connection 
$$D: C^\infty (X, E)\to C^\infty_{1} (X, E)$$ 
which can be extended as a linear differential operator of order 1 on the space of 
$(p,q)$ forms with value in $E$. Maybe the most familiar characterization 
of it is as follows:

\item {$\bullet $} We have the decomposition $D= \partial _h+ \dbar$ in types $(1,0)$ and respectively $(0,1)$; the $(0,1)$ part of $D$ is the extension of the usual $\dbar$ operator;
\smallskip
\item {$\bullet \bullet$} $D$ is compatible with the metric $h$.
\vskip 15pt
Unlike the operators $\dbar$, $\partial$ acting on $(p, q)$ forms, the components of $D$ do not commute in general. We measure the non commutativity defect by the {\it curvature} tensor, as follows. 

We denote by $[A, B]:= AB-(-1)^{\deg A \deg B}BA$ the graded commutator bracket of operators $A, B$; then it turns out that given $(E, h)$ a hermitian holomorphic vector bundle
the commutator
$$\Theta_h(E):= [\partial _h, \dbar]\in C^\infty(X, \Lambda^{1, 1}T^\star_X\otimes \Hom(E, E))$$
is of order zero and it is called the curvature tensor of $(E, h)$.

We define the operator $L_\omega: C^\infty(X, \Lambda^{p, q}T^\star_X\otimes E)\to
 C^\infty(X, \Lambda^{p+1, q+1}T^\star_X\otimes E)$ as follows
 $$Lu:= \omega\wedge u$$
 and its adjoint is denoted by $\Lambda_\omega$. It is easy to verify that with respect to some local trivialization of $E$, the expression of the adjoint operator is
 $$(\Lambda_\omega u)_{I \ol J}^\alpha= -\sqrt {-1}\sum_{i, j}\omega^{\ol j i}u^\alpha _{i\ol j I \ol J}$$
 where $\vert I\vert= p-1, \vert J\vert = q-1$ and the components $(u_{ K\ol L}^\alpha)$ of $u$ are assumed to be skew-symetric with respect to the ordered multi-index $K, L$.

\noindent Now the operators $\partial _h$ and $\dbar$ can be viewed as closed and densely defined operators on the Hilbert space $L^2(X, \Lambda ^{p, q}T^\star _X\otimes E)$. Their {\sl formal} adjoints 
$\displaystyle \partial _h^{\star } $ and $\dbar ^\star $ are densely defined and they admit extensions in the sense of distributions. Remark that in general they differ from the Hilbert space adjoints in the sense of 
von Neumann, as simple examples show it. However, it is well known that if the K\"ahler metric 
$\omega$ is geodesically complete, then the 2 operators coincide. We will use this fundamental fact in a moment, in the context of the pseudoconvex domains.

\noindent  Next we recall the fundamental Bochner-Kodaira-Nakano identity  ([4]); the context in which we will use it is the following.

We are given a K\"ahler weakly pseudoconvex manifold $(X, \omega)$ and a holomorphic line bundle $(L, h)$. Then for any  $(n, 1)$ form with compact support $\displaystyle u \in \cD(X, \Lambda ^{n, 1}T^\star _X\otimes L)$ we have the next relation
$$ \int_X \vert \dbar u\vert ^2_hdV_\omega+ \int_X \vert \dbar^\star u\vert ^2_hdV_\omega=
\int_X \vert \partial_h^\star u\vert ^2_hdV_\omega+ \int_X\langle [\Theta_h(L), \Lambda_\omega]u, u
\rangle _h dV_\omega $$
Since this relation holds true for any metric $h$ on $L$, we analyze next how does it changes when we consider a new metric $h_1:= \eta h$, where $\eta> 0$ is a smooth strictly positive function on $X$. We have
$$\eqalign{
\int_X \vert \dbar u\vert ^2_h\eta dV_\omega+ \int_X \vert \dbar^{\star_\eta} u\vert ^2_h \eta dV_\omega= & 
\int_X \vert \partial_\eta^\star u\vert ^2_h\eta dV_\omega+ \cr
+ & \int_X\langle [\Theta_h(L)-\sqrt {-1}\ddbar \log \eta, \Lambda_\omega]u, u
\rangle _h \eta dV_\omega\cr 
}$$

\noindent We use the identity
$$\ddbar \log \eta= {{\ddbar \eta}\over {\eta}}- {{\partial \eta\wedge \dbar \eta}\over {\eta ^2}}$$
and then the curvature term above becomes
$$  \int_X\langle [\eta \Theta_h(L)-\sqrt {-1}\ddbar \eta +
\sqrt {-1}{{\partial \eta\wedge \dbar \eta}\over {\eta }}, \Lambda_\omega]u, u
\rangle _h dV_\omega$$

The term $\displaystyle \int_X \vert \partial_\eta^\star u\vert ^2_h\eta dV_\omega$ is positive, therefore good enough for our future purposes; let us expand the term $\displaystyle \int_X \vert \dbar^{\star_\eta} u\vert ^2_h \eta dV_\omega$
(the reason is that later on we will choose $\eta$ converging to some function with logarithmic poles and we have to control in a very precise manner what will happen with this term under the convergence procedure). For any $u\in Dom(\dbar ^\star)$ and $v$ with compact support we have

$$\eqalign{\int_X\langle \dbar ^{\star_\eta}u, v\rangle _h\eta dV_\omega= & \int_X\langle u, \dbar v\rangle _h\eta dV_\omega= \cr
= & \int_X\langle u, \dbar (\eta v)\rangle _h dV_\omega- \int_X\langle u, \dbar \eta\wedge v\rangle _h dV_\omega= \cr
= &  \int_X\langle \dbar ^{\star}u, v\rangle _h\eta dV_\omega- \int_X\langle (\dbar \eta)^{\star}u, v\rangle _hdV_\omega\cr }.$$

In conclusion, the formal adjoint operator corresponding to the twisted metric $h_1$ is equal to
$$\dbar ^{\star_\eta}u= \dbar ^{\star}u- 1/\eta  (\dbar \eta)^{\star}u$$
and let us now compute the corresponding term in the Bochner formula
$$\eqalign{ 
\int_X \vert \dbar^{\star_\eta} u\vert ^2_h \eta dV_\omega= & 
\int_X \vert \dbar ^{\star}u- 1/\eta  (\dbar \eta)^{\star}u \vert ^2_h \eta dV_\omega= \cr
= & \int_X \vert \dbar^{\star} u\vert ^2_h \eta dV_\omega + 
\int_X {1\over \eta} \vert (\dbar \eta)^{\star}u \vert ^2_h dV_\omega- 2\Re 
\int_X\langle \dbar ^{\star}u, (\dbar \eta)^{\star}u\rangle _hdV_\omega\cr
}$$

We claim next that we have the following identity
$$\vert (\dbar \eta)^{\star}u \vert ^2_h=  \langle [\sqrt {-1}\partial \eta\wedge \dbar \eta,  \Lambda_\omega]u, u\rangle _h\leqno (7)$$
at each point of $X$. Indeed, let us check it by a computation in local coordinates.
We take $x\in X$ an arbitrary point and let $(z^j)$ be local coordinates on some open set $\Omega$ centered at $x$, which are geodesic for the K\"ahler metric $\omega$ at this point. Locally we have
$$u= \sum_j u_{\ol j}dz\wedge dz^{\ol j}\otimes e_\Omega$$
where $dz:= dz^1\wedge...\wedge dz^n$ and $e_\Omega$ is a local holomorphic frame of $L$.
Then by the definition of $(\dbar \eta)^\star$ we have
$$(\dbar \eta)^\star(u)= (-1)^n\sum_ju_{\ol j}{{\partial \eta}\over {\partial z^j}}dz\otimes e_\Omega$$
at $x$, and therefore we get 
$$\vert (\dbar \eta)^{\star}u \vert ^2_h= \sum_{j,k}u_{\ol j}\ol {u_{\ol k}}{{\partial \eta}\over {\partial z^j}}{{\partial \eta}\over {\partial z^{\ol k}}}\exp\big(-\varphi_L(x)\big)\leqno (8)$$
On the other hand, the expression of  the contraction operator $\Lambda_\omega$ on $u$ reads as
$$\Lambda_\omega(u)= \sqrt {-1}\sum _j(-1)^{n+j+1}u_{\ol j}dz^1\wedge \wh {dz^j}\wedge dz^n
\otimes e_\Omega$$
and we have
$$\eqalign{
\sqrt {-1}\partial \eta\wedge \dbar \eta\wedge  \Lambda_\omega (u)= & \sum_{j, p, q}(-1)^{n+ j}
u_{\ol j}{{\partial \eta}\over {\partial z^p}}{{\partial \eta}\over {\partial z^{\ol q}}}
dz^p\wedge dz^{\ol q}\wedge dz^1\wedge \wh {dz^j}\wedge dz^n
\otimes e_\Omega = \cr 
= & \sum_{j, k}u_{\ol j}{{\partial \eta}\over {\partial z^j}}{{\partial \eta}\over {\partial z^{\ol k}}}
 dz\wedge dz^{\ol k}
\otimes e_\Omega.\cr
}$$
The relation (7) follows, as by the previous identity we have
$$\langle [\sqrt {-1}\partial \eta\wedge \dbar \eta,  \Lambda_\omega]u, u\rangle _h= \sum_{j,k}u_{\ol j}\ol {u_{\ol k}}{{\partial \eta}\over {\partial z^j}}{{\partial \eta}\over {\partial z^{\ol k}}}\exp\big(-\varphi_L(x)\big).\leqno (9)$$

Summing up what we have obtained so far, the twisted Bochner identity become
$$\eqalign{ 
\int_X \vert \dbar u\vert ^2_h\eta dV_\omega+ \int_X \vert \dbar^{\star} u\vert ^2_h \eta dV_\omega= & 
\int_X \vert \partial_\eta^\star u\vert ^2_h\eta dV_\omega+ \cr
+ & \int_X\langle [\eta \Theta_h(L)-\sqrt {-1}\ddbar \eta, \Lambda_\omega]u, u
\rangle _h dV_\omega+\cr 
+ & 2\Re 
\int_X\langle \dbar ^{\star}u, (\dbar \eta)^{\star}u\rangle _hdV_\omega
\cr 
}$$

Let $\lambda$ be another positive function on $X$; by completing the square in the above identity and neglecting the positive term in the right hand side term above, we have proved the next inequality (see e.g. [7], [8], [10], [12]).

\claim Lemma|Let $\eta, \lambda$ be smooth functions on $X$, and let $\displaystyle u \in \cD(X, \Lambda ^{n, 1}T^\star _X\otimes L)$. Then  we have the next inequality 

$$\eqalign{ 
\int_X \vert \dbar u\vert ^2_h \eta dV_\omega+ & \int_X (\eta +\lambda)\vert \dbar^{\star} u\vert ^2_h dV_\omega\geq  \cr  
\geq & \int_X\langle \big[\eta \Theta_h(L)-\sqrt {-1}\ddbar \eta-
\sqrt {-1}{{\partial \eta\wedge \dbar \eta}\over {\lambda }}, \Lambda_\omega\big]u, u
\rangle _h dV_\omega\cr 
}
\leqno (10)$$
\endclaim 

\section{\S D Functional analysis}

\vskip 10pt
\noindent We consider the modified $\dbar$ operators 
$$Tu:= \dbar \big(\sqrt {\eta + \lambda}u\big) \hbox { and } Su:= \sqrt \eta\big(\dbar u\big)$$
acting on $(p, q)$ forms with values in a bundle $F$. They are densely defined, and we obviously have
$S\circ T= 0$.

Then the previous lemma can be reformulated as follows: {\sl for any $(n, 1)$--form $u$
with values in a line bundle $(F, h)$ which belongs to the domains of $S$ and $T^\star$ we have}
$$\Vert T^\star u\Vert^2+ \Vert Su\Vert ^2\geq 
\langle\! \langle \big[\eta \Theta_h(F)-\sqrt {-1}\ddbar \eta-
{{\partial \eta\wedge \dbar \eta}\over {\lambda }}, \Lambda_\omega\big]u, u
\rangle\!\rangle. \leqno (11)$$
We stress on the fact that this was only proved for compactly supported forms $u$; in general, we use a standard 
density argument--remark that it is at this point that we need the 
fact that $X$ carries a complete K\"ahler manifold (see [2]).

Now assume that for some specific choice of the functions $\eta$ and $\lambda$ above we have
$$\eta \Theta_h(L)-\sqrt {-1}\ddbar \eta-
{{\partial \eta\wedge \dbar \eta}\over {\lambda }}\geq \sqrt {-1}\tau\partial \mu\wedge \dbar \mu- \delta\omega\leqno (12)$$
on $X$; here $\tau$ is a positive function, $\mu$ is an arbitrary one,
and $\delta> 0$ is a positive real number. Then the curvature term in the twisted Bochner inequality above verify the next relation

$$\langle\! \langle \big[\eta \Theta_h(L)-\sqrt {-1}\ddbar \eta-
{{\partial \eta\wedge \dbar \eta}\over {\lambda }}, \Lambda_\omega\big]u, u
\rangle\!\rangle\geq \Vert \sqrt \tau (\dbar \mu)^\star u\Vert ^2
-\delta \Vert u\Vert^2
$$
and therefore the relation (11) become 
$$\Vert T^\star u\Vert^2+ \Vert Su\Vert ^2+ \delta \Vert u\Vert^2\geq \Vert \sqrt \tau (\dbar \mu)^\star u\Vert ^2
\leqno (13)$$

\noindent  As usual, we want to solve the equation $Tv=g$ for some closed $(n, 1)$ form $g$ with values in $F:= Z+ L$, such that $Sg= 0$. Given the lower bound $(13)$ above (the corresponding operator is {\sl never} positively defined,
even if $\delta= 0$), we cannot hope to solve the previous equation for all $g$ in the kernel of $\dbar$, but nevertheless let us consider 
$g= \dbar \tau\wedge g_0+ g_1$, where $g_0$ is a $L^2$ form of $(n, 0)$ type,
and $g_1$ is a $(n, 1)$-form.

\noindent Then we have
$$\eqalign{
{1\over 2}\Big\vert \int_X\langle g, u\rangle _hdV_\omega\Big\vert ^2\leq  & \Big\vert \int_X\langle  \dbar \tau\wedge g_0, u\rangle _hdV_\omega\Big\vert ^2
+ \Big\vert \int_X\langle g_1, u\rangle _hdV_\omega\Big\vert ^2= \cr
= & \Big\vert \int_X\langle  g_0, (\dbar \tau)^\star u\rangle _hdV_\omega\Big\vert ^2+  \Big\vert \int_X\langle g_1, u\rangle _hdV_\omega\Big\vert ^2\leq \cr
 \leq & \Vert \sqrt \tau (\dbar \tau)^\star u\Vert ^2\int_X 1/\tau \vert g_0\vert ^2 dV_\omega+ \delta\Vert u\Vert^2\cdot {1\over \delta}\Vert g_1\Vert^2\leq \cr
 \leq & C(g, \tau, \delta)\Big(\Vert T^\star u\Vert^2+ \delta \Vert u\Vert^2\Big);\cr   
}
\leqno (14)$$
here we use the standard trick --see [1], [6]-- in order to avoid the term $\Vert Su\Vert^2$ in the above inequality. We have employed the notation 
$$C(g, \tau,\delta):= \int_X 1/\tau \vert g_0\vert ^2 e^{-\varphi_F}+ {1\over \delta}
\int_X|g_1|^2_\omega e^{-\varphi_F} dV_\omega.$$
\vskip 7pt
\noindent In conclusion, under the precise circumstances described in this paragraph the map
$$(T^\star u, \sqrt \delta u)\to \langle g, u\rangle$$
is continuous, hence 
we can solve the approximate 
$\dbar$ equation 
$$Tv+ \sqrt \delta w= \dbar \tau\wedge g_0+ g_1$$ 
together with the estimates
$$\int_X|v|^2e^{-\varphi_F}+ \int_X|w|^2_\omega e^{-\varphi_F}dV_\omega\leq C(g, \tau, \delta).\leqno (15)$$
We remark that all the computations/considerations above 
were done under the assumption that the metric 
$\omega$ is complete. 

\section {\S E End of the proof}
\vskip 10pt

Since we want the final estimates of the extension to depend on the restriction of $U_\infty$ to $Z$ only, it is natural to truncate the section $U_\infty$ as follows: let $\theta: \bR\to [0, 1]$ be a smooth function with support in $]-\infty, 1[$, such that $\theta\equiv 1$ on the interval 
$]-\infty, 1/2]$, and such that $\sup_\bR|\theta^\prime|\leq 4$. Then we define 
$$g_\varepsilon:= \dbar \Big(\theta\big(\vert s\vert^2/\varepsilon ^2\big)U_\infty\Big)$$
It is an $(n, 1)$ form with values in $F:= E+ L$, which is $L^2$ with respect to the original metric 
$h_L$ of the bundle $L$, twisted with the singular weight function $\log|s|^2$ (we will verify this assertion in a moment). Let $\chi_0$ be any function on $\bR_-$
whose derivative is strictly positive in absolute value (we will make a specific choice later). We compute
$$\eqalign{
g_\varepsilon = & {{1}\over {\varepsilon ^2}}\theta^\prime \Big({{\vert s\vert^2}\over {\varepsilon ^2}}\Big)\dbar |s|^2\wedge U_\infty + 
\theta\big(\vert s\vert^2/\varepsilon ^2\big)\dbar U_\infty
=\cr
= &\big(1+ {{|s|^2}\over {\varepsilon ^2}}\big){{\theta^\prime\Big({{\vert s\vert^2}\over {\varepsilon ^2}}\Big)
}\over {\chi_0^\prime\big(\log(\varepsilon ^2+ \vert s\vert^2)\big)}}\dbar \Big(\chi_0\big(\log(\varepsilon ^2+ \vert s\vert^2)\big)\Big)\wedge U_\infty+  
\theta\big(\vert s\vert^2/\varepsilon ^2\big)\dbar U_\infty:= \cr
= & \dbar \eta_\varepsilon\wedge g_{\varepsilon, 0}+ g_{\varepsilon, 1}\cr
}$$
 
The motivation for introducing the logarithmic term above is that in this way the function
in front of the $\dbar$ is bounded; also, it is coherent with the "functional analysis" discussion above. In the last line of the above equality we have used the function 
$$\eta_\varepsilon: = -\chi_0\big(\log(\varepsilon ^2+ \vert s\vert^2)\big)
\leqno (16)$$ and the $(n, 0)$--form  $g_{\varepsilon, 0}$ which is essentially $U_\infty$ up to a bounded function.  

According to the previous section D (see especially the relations (12), (13)), the lower bound we have to obtain for the curvature term in (12) must 
be something of the form $\tau\sqrt {-1}\partial \eta_\varepsilon \wedge \dbar \eta_\varepsilon$, so we are going to expand the term 
$$\eta_\varepsilon(\Theta_h(F)+ \sqrt {-1}\ddbar \log|s|^2)-\sqrt {-1}\ddbar \eta_\varepsilon$$
and choose an appropriate function $\lambda_\varepsilon$ in order to obtain the 
correct lower bound of the curvature.

We define the function
$$\sigma_\varepsilon:=  \log(\varepsilon ^2+ \vert s\vert^2);$$ 
then a straightforward computation shows that 
$$\sqrt{-1}\ddbar\sigma_\varepsilon \geq \sqrt {-1}{{\varepsilon ^2}\over {|s|^2}}
\partial \sigma_\varepsilon \wedge \dbar \sigma_\varepsilon-
{{\langle \Theta(E)s, s\rangle}\over {\varepsilon^2+ |s|^2}}\leqno (17)$$ 
and we clearly have 
$$-\ddbar \eta_\varepsilon= \chi_0^\prime(\sigma_\varepsilon) \ddbar \sigma_\varepsilon+ 
\chi_0^{\prime\prime}(\sigma_\varepsilon)\partial \sigma_\varepsilon\wedge \dbar \sigma_\varepsilon.\leqno(18)$$

We want to use the relation (17) in order to get a lower bound for the hessian of $-\eta_\varepsilon$, so we need some definite positive 
bounds for the derivative of $\chi_0$; we assume $1\leq \chi_0^\prime\leq 2$. Then we have  
$$-\sqrt{-1}\ddbar \eta_\varepsilon\geq \Big( {{\varepsilon ^2}\over {2|s|^2}}+ {{\chi_0^{\prime\prime}(\sigma_\varepsilon)}\over {\chi_0^\prime(\sigma_\varepsilon)^2}}
 \Big)\sqrt{-1}\partial\eta_\varepsilon\wedge \dbar \eta_\varepsilon- {{\chi_0^\prime(\sigma_\varepsilon)}\over {\varepsilon^2+ |s|^2}}
 \langle \Theta(E)s, s\rangle \leqno (19)$$
\medskip

We use now the regularized metric $\varphi_{L, \rho}$ of $L$ obtained in the paragraph A; we recall here its curvature properties

{\itemindent 3mm 
\smallskip
\item {$\bullet$} $\displaystyle \Theta_{h_{L, \rho}}(L)\geq  -\delta_\rho\omega$;
\smallskip
\item {$\bullet$} $\displaystyle \Theta_{h_{L, \rho}}(L)\geq  
{1\over \alpha}
\Theta (E)- \delta_\rho\omega$ 
for some $\alpha\geq
1$.
}

\noindent By the normalization condition of the section $s$ we have 
$\eta_\varepsilon\geq 2\alpha$; in particular
we get $\displaystyle \eta_\varepsilon\geq  {{\alpha\chi_0^\prime (\sigma_\varepsilon)|s|^2}\over {\varepsilon^2+ |s|^2}}$
and therefore we get the inequality 
$$
\eqalign{
\eta_\varepsilon & (\Theta_{h_{\rho}}(L)+ \Theta_{h}(E)+ 
\sqrt {-1}\ddbar \log|s|^2))\geq \cr
\geq & {{\chi_0^\prime(\sigma_\varepsilon)}\over {\varepsilon^2+ |s|^2}}
\Theta(E)- \delta_\rho\eta_\varepsilon \omega
}\leqno (20) $$
(remark that we have used at this point both curvature assumptions above).

By the relation (18) combined with (20), it is clear that the ideal candidate for the function $\lambda_\varepsilon$
would be 
$$\lambda_\varepsilon:= {{\chi_0^\prime(\sigma_\varepsilon)^2}\over {\chi_0^{\prime\prime}(\sigma_\varepsilon)}}$$
provided that the denominator is non-zero. A choice of function $\chi_0$ which will be convenient for our purposes is as in [4]
$$\chi_0(t) := t-\log(1-t)\leqno (21)$$
(we refer to the article [9] for a very complete treatment of other possible choices of $\chi_0$).

By collecting all the relations above we get
$$\eqalign{
\eta_\varepsilon & (\Theta_{h_{\rho}}(L)+ \Theta_{h}(E)+ 
\sqrt {-1}\ddbar \log|s|^2)-\cr
- & \sqrt {-1}\ddbar \eta_\varepsilon-{{\sqrt {-1}}\over {\lambda_\varepsilon}}\partial \eta_\varepsilon\wedge 
\dbar \eta_\varepsilon
\geq {{\varepsilon ^2}\over {2|s|^2}}\sqrt {-1}\partial \eta_\varepsilon\wedge 
\dbar \eta_\varepsilon - \delta_\rho\eta_\varepsilon \omega. \cr
}
\leqno (22)$$

\noindent The above inequality holds on $\displaystyle X_\rho:= X\setminus 
\Big(E_\rho \big(\Theta_{h_L}(L)\big)\cup Z\Big)$, and it involves the metric $\omega$ which is not complete on $X_\rho$. We remark that the completeness of the metric is necessary in order to use the results obtained in the functional analysis 
paragraph D, hence we will use next the following lemma (cf. [2]).
\medskip

\claim Lemma {[2]}|Let $(X, \omega)$ be a compact K\"ahler manifold, and let $W\subset X$ be an 
analytic subset. Then $X\setminus W$ carries a complete K\"ahler metric.
\endclaim

\medskip

\noindent Let $\omega_\rho$ be a complete metric on $X_\rho$; for each $k\geq 1$, we define the 
metric 
$$\omega_k:= \omega+ {1\over k}\omega_\rho \leqno (23)$$
it is equally complete on $X_\rho$, and we have $\omega_k> \omega$.
Therefore the inequality (22) is still valid if we are using $\omega_k$
instead of $\omega$, i.e. we have

$$\eqalign{
\eta_\varepsilon & (\Theta_{h_{\rho}}(L)+ \Theta_{h}(E)+ 
\sqrt {-1}\ddbar \log|s|^2)-\cr
- & \sqrt {-1}\ddbar \eta_\varepsilon-{{\sqrt {-1}}\over {\lambda_\varepsilon}}\partial \eta_\varepsilon\wedge 
\dbar \eta_\varepsilon
\geq {{\varepsilon ^2}\over {2|s|^2}}\sqrt {-1}\partial \eta_\varepsilon\wedge 
\dbar \eta_\varepsilon - \delta_\rho\eta_\varepsilon \omega_k. \cr
}
\leqno (24)$$
for each $k\geq 1$. 

For each $\varepsilon> 0$ we can choose $\rho_\varepsilon> 0$ 
such that 
$$\delta_{\rho_\varepsilon}\log^2{1\over \varepsilon^2}\leq \varepsilon^{\varepsilon_0\over 1+ \varepsilon_0}$$
(see the inequality (6) in the paragraph B) and then the inequality (24) implies
$$\eqalign{
\eta_\varepsilon & (\Theta_{h_{\rho_\varepsilon}}(L)+ \Theta_{h}(E)+ 
\sqrt {-1}\ddbar \log|s|^2)-\cr
- & \sqrt {-1}\ddbar \eta_\varepsilon-{{\sqrt {-1}}\over {\lambda_\varepsilon}}\partial \eta_\varepsilon\wedge 
\dbar \eta_\varepsilon
\geq {{\varepsilon ^2}\over {2|s|^2}}\sqrt {-1}\partial \eta_\varepsilon\wedge 
\dbar \eta_\varepsilon - \varepsilon^{\varepsilon_0\over 1+ \varepsilon_0} \omega_k, \cr
}
\leqno (25)$$
given the expression of $\eta_\varepsilon$. 

\smallskip

We see that all the conditions in the ``Functional analysis" paragraph D are fulfilled and 
as a consequence, for each $\varepsilon> 0$ we can solve the following approximate $\dbar$ equation
$$\dbar \big(\sqrt {\eta_\varepsilon+ \lambda_\varepsilon}v_{\varepsilon, k}\big)+ \varepsilon^{\varepsilon_0\over 2(1+ \varepsilon_0)} w_{\varepsilon, k}= \dbar \eta_\varepsilon\wedge g_{\varepsilon, 0}+ g_{\varepsilon, 1}; $$
in addition, we have the estimates
$$\eqalign{
\int_X& {{|v_{\varepsilon, k}|^2}\over 
{|s|^{2}}}e^{-\varphi_E- \varphi_{L, \rho_\varepsilon}}+
\int_X{{|w_{\varepsilon, k}|_{\omega_k}^2}\over 
{|s|^{2}}}e^{-\varphi_E- \varphi_{L, \rho_\varepsilon}}dV_{\omega_k}
\leq \cr
\leq & C\int_X\theta^\prime \Big({{\vert s\vert^2}\over {\varepsilon ^2}}\Big)^2{{|U_\infty|^2}\over {|s|^2}}
e^{-\varphi_{L, \rho_\varepsilon}-\varphi_E}+ \cr
+ & \varepsilon^{-{\varepsilon_0\over 1+ \varepsilon_0}}
\int_X \theta\Big({{\vert s\vert^2}\over {\varepsilon ^2}}\Big)^2
{|\dbar U_\infty|^2_{\omega_k}\over 
\vert s\vert^2}e^{-\varphi_{L, \rho_\varepsilon}-\varphi_E}dV_{\omega_k}
\cr
}
\leqno (26)$$
\medskip

\noindent We analyze next the limit as $\varepsilon\to 0$ of the right hand side
terms of the inequality (26). By the properties of the regularization family 
$(\varphi_{L, \rho})$, we obtain the relation
$$\int_X\theta^\prime \Big({{\vert s\vert^2}\over {\varepsilon ^2}}\Big)^2{{|U_\infty|^2}\over {|s|^2}}
e^{-\varphi_{L, \rho_\varepsilon}-\varphi_E}\leq 
\int_X\theta^\prime \Big({{\vert s\vert^2}\over {\varepsilon ^2}}\Big)^2{{|U_\infty|^2}\over {|s|^2}}
e^{-\varphi_{L}-\varphi_E}.\leqno (27)$$
At this point we have to use the hypothesis $\cC$ concerning the singularities of $\varphi_L$: it implies that the limit as 
$\varepsilon\to 0$ of the right hand side term of (27) is equal to
$$\int_Z|u|^2e^{-\varphi_L}\leqno (28)$$
up to a universal constant.

As for the remaining term, we observe that we have
$$\int_X \theta\Big({{\vert s\vert^2}\over {\varepsilon ^2}}\Big)^2
{|\dbar U_\infty|^2_{\omega_k}\over 
\vert s\vert^2}e^{-\varphi_{L, \rho_\varepsilon}-\varphi_E}dV_{\omega_k}\leq 
\int_X \theta\Big({{\vert s\vert^2}\over {\varepsilon ^2}}\Big)^2
{|\dbar U_\infty|^2_{\omega}\over 
\vert s\vert^2}e^{-\varphi_{L, \rho_\varepsilon}-\varphi_E}dV_{\omega}$$
since $\omega_k> \omega$ and in this context we have the following result,
due to [2].

\claim Lemma([2])|
Let $\omega_{1}\leq\omega_{2}$ be two K\"{a}hler metrics and let
$\|\cdot\|_{1,2}$ be the corresponding norms. Then, if $f$ is an
$(n,q)$-form we have the pointwise inequality
$$\|f\|_{2}^{2}dV_{\omega_{2}}\leq\|f\|_{1}^{2}dV_{\omega_{1}}.$$
\endclaim

\noindent By the arguments we have already used a few lines ago, we obtain

$$\int_X \theta\Big({{\vert s\vert^2}\over {\varepsilon ^2}}\Big)^2
{|\dbar U_\infty|^2_{\omega}\over 
\vert s\vert^2}e^{-\varphi_{L, \rho_\varepsilon}-\varphi_E}dV_{\omega}\leq \int_X \theta\Big({{\vert s\vert^2}\over {\varepsilon ^2}}\Big)^2
{|\dbar U_\infty|^2_{\omega}\over 
\vert s\vert^2}e^{-\varphi_{L}-\varphi_E}dV_{\omega},$$
and the relation (6) of the paragraph B shows that we have
$$
\int_X \theta\Big({{\vert s\vert^2}\over {\varepsilon ^2}}\Big)^2
{|\dbar U_\infty|^2_{\omega}\over |s|^2}e^{-\varphi_{L}-\varphi_E}dV_{\omega}\leq
C\varepsilon^{2\varepsilon_0\over 1+ \varepsilon_0}\log ^{2\over 1+ \varepsilon_0}{1\over \varepsilon^2}\leqno (29)$$

\noindent Combining the inequalities (26), (27) and (29) we obtain
$$\eqalign{
\int_X& {{|v_{\varepsilon, k}|^2}\over 
{|s|^{2}}}e^{-\varphi_E- \varphi_{L, \rho_\varepsilon}}+
\int_X{{|w_{\varepsilon, k}|_{\omega_k}^2}\over 
{|s|^{2}}}e^{-\varphi_E- \varphi_{L, \rho_\varepsilon}}dV_{\omega_k}
\leq \cr
\leq & C_0\int_Z|u|^2e^{-\varphi_L}+ C_\varepsilon\cr
}
\leqno (30)$$
where $C_\varepsilon\to 0$ as $\varepsilon\to 0$ and $C_0$ is a numerical constant.
\medskip

\noindent We intend to remove the dependence in $k$ from the relation (30), by extracting a limit, as follows (the procedure is completely similar to the one employed in [2]).
 
For each $\varepsilon> 0$ fixed, the $(n, 0)$-forms $(v_{\varepsilon, k})$
are defined on the manifold $X\setminus (Z\cup Y_\varepsilon)$ 
and bounded in $L^2$ norm with respect to the metric $\varphi_{L, \rho_\varepsilon}+ \log |s|^2$ by a constant which is independent of $k$. Therefore we can extract a weak limit say 
$$v_\varepsilon\in L^2\big(X\setminus (Z\cup Y_\varepsilon)\big)$$ 
from the sequence
$(v_{\varepsilon, k})$.  

The arguments we use in order to derive a similar conclusion 
for the sequence of $(n, 1)$--forms 
$(w_{\varepsilon, k})$ are slightly more technical, as follows.
Let $k\geq l$ be two integers; since $\omega_l\geq \omega_k$, we have
$$\int_X{{|w_{\varepsilon, k}|_{\omega_l}^2}\over 
{|s|^{2}}}e^{-\varphi_E- \varphi_{L, \rho_\varepsilon}}dV_{\omega_l}
\leq \int_X{{|w_{\varepsilon, k}|_{\omega_k}^2}\over 
{|s|^{2}}}e^{-\varphi_E- \varphi_{L, \rho_\varepsilon}}dV_{\omega_k}
\leq C\leqno(31)$$
as one can infer thanks to (30). 
Therefore, we can extract a weak limit $w_\varepsilon$ of $(w_{\varepsilon, k})$; for any compact $K\subset X\setminus (Z\cup Y_\varepsilon)$, the metric 
$\omega_l$ is comparable with the non-complete metric $\omega$, so
in conclusion, we let $k\to \infty$ in (30) and we obtain the estimate
$$\eqalign{
\int_X& {{|v_{\varepsilon}|^2}\over 
{|s|^{2}}}e^{-\varphi_E- \varphi_{L, \rho_\varepsilon}}+
\int_X{{|w_{\varepsilon}|_{\omega}^2}\over 
{|s|^{2}}}e^{-\varphi_E- \varphi_{L, \rho_\varepsilon}}dV_{\omega}
\leq \cr
\leq & C_0\int_Z|u|^2e^{-\varphi_L}+ C_\varepsilon\cr
}
\leqno (32)$$
together with the equality
$$\dbar \big(\sqrt {\eta_\varepsilon+ \lambda_\varepsilon}v_{\varepsilon}\big)+ \varepsilon^{\varepsilon_0\over 2(1+ \varepsilon_0)} w_{\varepsilon}= 
\dbar \Big(\theta\Big({|s|^2\over \varepsilon^2}\Big)U_\infty\Big)\leqno(33)$$
in the sense of distributions on $X\setminus (Z\cup Y_\varepsilon)$.
The following lemma (cf. [2] and the references therein) shows that the equation (33) is verified on the 
manifold $X$ in the sense of distributions.

\claim Lemma ([2])|
Let $\Omega$ be an open subset of $\bC^{n}$ and let $Y$ be a complex
analytic subset of $\Omega$. Assume that $v$ is a (p,q-1)-form with
$L_{loc}^{2}$ coefficients and $w$ a (p,q)-form with $L_{loc}^{1}$
coefficients such that $\overline{\partial}v=w$ on $\Omega\backslash
Y$ (in the sense of distributions). Then $\overline{\partial}v=w$ on
$\Omega$.
\endclaim
\medskip
Besides the fact that $v_\varepsilon$ and $w_\varepsilon$ are in 
$L^2$, we do not have further informations about their regularity. Indeed, 
from the estimate (32) one would like to infer that $v_\varepsilon$ is equal to zero
when restricted to $Z$, but in the $L^2$ setting, this makes no sense.
We tackle this difficulty as follows: the equality (33) shows that we have
$$\dbar w_\varepsilon= 0$$ 
so locally on each coordinate polydisk $\Omega\subset X$ we can
solve the equation
$$\dbar f_{\Omega, \varepsilon}= w_\varepsilon$$
such that we have
$$\int_\Omega {|f_{\Omega, \varepsilon}|^2\over |s|^2}e^{-\varphi_E- \varphi_{L, \rho_\varepsilon}}d\lambda\leq C
\int_X{{|w_{\varepsilon}|_{\omega}^2}\over 
{|s|^{2}}}e^{-\varphi_E- \varphi_{L, \rho_\varepsilon}}dV_{\omega}\leqno (34)$$
by H\"ormander $L^2$ estimates (cf. [2], [6]). We consider a finite covering
of $X$ with polydisks $\Omega_j$, and we denote by $f_{j, \varepsilon}$
the solutions of the corresponding equation (34).
Then we have
$$\dbar\Big(\theta\Big({|s|^2\over \varepsilon^2}\Big)U_\infty- 
\sqrt{\eta_\varepsilon+ \lambda_\varepsilon}v_\varepsilon
-\varepsilon^{\varepsilon_0\over 2(1+ \varepsilon_0)} f_{j, \varepsilon}\Big)= 0$$
in other words, the local section 
$$U_{j, \varepsilon}:= \theta\Big({|s|^2\over \varepsilon^2}\Big)U_\infty- 
\sqrt{\eta_\varepsilon+ \lambda_\varepsilon}v_\varepsilon
-\varepsilon^{\varepsilon_0\over 2(1+ \varepsilon_0)} f_{j, \varepsilon}$$
of the bundle $K_X+ Z+ L$ defined on $\Omega_j$ is holomorphic. In particular, this implies that the function
$$\sqrt{\eta_\varepsilon+ \lambda_\varepsilon}v_\varepsilon
+\varepsilon^{\varepsilon_0\over 2(1+ \varepsilon_0)} f_{j, \varepsilon}$$
is $\cC^\infty(\Omega_j)$, and the estimates (32) and (34) show that this function 
is equal to zero when restricted to $Z\cap \Omega_j$. Hence we have
$$U_{j, \varepsilon}|_{Z\cap \Omega_j}= u$$
and we will show next that:

\noindent $\bullet$ As $\varepsilon\to 0$, the local sections piece together, and 
give an extension of $u$.

\noindent $\bullet$ The extension obtained as described in the previous bullet
verifies the estimates required by the main theorem.

\medskip 
The verification is quite easy: in the expression of $U_{j, \varepsilon}$, the only 
``local" component is the term 
$$\varepsilon^{\varepsilon_0\over 2(1+ \varepsilon_0)} f_{j, \varepsilon};$$
by the estimate (34), there exists $f_j\in L^2(\Omega_j, e^{-\varphi_L-\log|s|^2})$ such that $f_{j, \varepsilon}\to f_j$ weakly. But then we have
$$\varepsilon^{\varepsilon_0\over 2(1+ \varepsilon_0)} f_{j, \varepsilon}\to 0$$
as $\varepsilon\to 0$, and this shows that we will have
$$U_j|_{\Omega_j\cap \Omega_l}= U_l|_{\Omega_j\cap \Omega_l}$$
for any pair of indexes $(j, l)$, where $U_j$ is a 
limit of $U_{j, \varepsilon}$. We denote by $U$ the resulting section.

In order to evaluate the $L^2$ norm of $U$, we remark that we have
 $$\eta_\varepsilon+ \lambda_\varepsilon\leq 5\log (|s|^{2}+ \varepsilon ^2)^2$$ for $\varepsilon \ll 1$; we obtain 
$$\eqalign{\int_{\Omega_j}{{|U_{j, \varepsilon}|^2  e^{-\varphi_{L, \rho_\varepsilon}-
\varphi_E}}\over {(|s|^{2}+ \varepsilon ^2)\big(\log (|s|^{2}+ \varepsilon ^2)\big)^2}}\leq &  C_0\int_Z |u|^2e^{-\varphi_L- \varphi_E}+ \cr
+ & C_\varepsilon+ \int_X{|\theta(|s|^2/\varepsilon^2)U_\infty|^2e^{-\varphi_{L, \rho_\varepsilon}-
\varphi_E}\over
(|s|^{2}+ \varepsilon ^2)\big(\log (|s|^{2}+ \varepsilon ^2)\big)^2}\cr
}
\leqno(35)$$
where $C_\varepsilon\to 0$ as $\varepsilon \to 0$ (cf. (34) combined with (30)).

The proof will be finished if we can show that the term
$$\int_X{|\theta(|s|^2/\varepsilon^2)U_\infty|^2e^{-\varphi_{L, \rho_\varepsilon}-
\varphi_E}\over
(|s|^{2}+ \varepsilon ^2)\big(\log (|s|^{2}+ \varepsilon ^2)\big)^2}$$
converges to zero. This is a consequence of the H\"older formula, as follows.
In the first place, given the expression of the section
$U_\infty$, we see that it would be enough to show that we have
$$\int_{\Omega_\alpha}{|\theta(|s|^2/\varepsilon^2)|\wt f_\alpha|^2e^{-\varphi_{L}-
\varphi_E}\over
(|s|^{2}+ \varepsilon ^2)\big(\log (|s|^{2}+ \varepsilon ^2)\big)^2}\to 0\leqno (36)$$
as $\varepsilon\to 0$.
We have
$$\eqalign{
\int_{\Omega_\alpha}{|\theta(|s|^2/\varepsilon^2)\wt f_\alpha|^2e^{-\varphi_{L}-
\varphi_E}\over
(|s|^{2}+ \varepsilon ^2)\big(\log (|s|^{2}+ \varepsilon ^2)\big)^2}\leq & 
\Big(C\int_{\Omega_\alpha}{|\wt f_\alpha|^{2}e^{-(1+ \varepsilon_0)\varphi_{L}}\over
|s|^{2}\big(\log (|s|^{2})\big)^2}\Big)^{1\over 1+\varepsilon_0}\cr
& \Big(\int_{\Omega_\alpha\cap (|s|< \varepsilon)}{d\lambda\over
(|s|^{2}+ \varepsilon ^2)\big(\log (|s|^{2}+ \varepsilon ^2)\big)^2}
\Big)^{\varepsilon_0\over 1+\varepsilon_0}\cr
}
\leqno (37)$$
and so we see that we have
the inequality
$$\int_X{|\theta(|s|^2/\varepsilon^2)U_\infty|^2e^{-\varphi_{L, \rho_\varepsilon}-
\varphi_E}\over
(|s|^{2}+ \varepsilon ^2)\big(\log (|s|^{2}+ \varepsilon ^2)\big)^2}\leq
{C\over \log^{2\varepsilon_0\over 1+\varepsilon_0}({1\over \varepsilon^2})}$$
and we let $\varepsilon\to 0$ in relation (35) in order to obtain precisely the estimate claimed in the main statement.
Therefore, we are done.\hfill \qed
\vskip 20pt

\section{\S F Comments about the general case}

\vskip 7pt
The extension result 0.1 is of course expected to be proved
in a more general setting:
it should be enough to assume that $X$ is a compact K\"ahler manifold.
If the metric $h_L= e^{-\varphi_L}$ has logarithmic poles, then our main theorem is enough to prove the result --this is just because in this case,
our additional hypothesis are easily checked. 

As we can see from the proof of Theorem 0.1, the hypothesis
$$u\in {\cal I}_+(h_L|_Z)$$
together with the hypothesis $\cal C$ are needed because we cannot allow
the parameters $\rho$ (appearing in the regularization process) and
$\varepsilon$ (whose effect is to concentrate the mass along $Z$ in the estimate (35)) to vary independently. We describe next the problems we encounter while trying to 
bypass this difficulty.

By using the regularization in [3] by applying just the convolution kernel
(i.e. without the Legendre transform), we obtain a family of non-singular
metrics $h_{L, \rho}$ on $L$ with the following properties

{\itemindent 0.5cm

\item{(1)} $\displaystyle \Theta_{h_{L, \rho}}(L)\geq -\mu_\rho\omega$;
\smallskip
\item{(2)} $\displaystyle \Theta _{h_{L, \rho}}(L) \geq
{1\over \alpha}\Theta(E) -\mu_\rho\omega$ for some
$\alpha\geq 1, \quad \alpha\in{\cC}^\infty(X)$.}

\noindent In the inequality above, $(\mu_\rho)_{\rho> 0}$ is a family of 
positive functions on $X$, such that:

\noindent $\bullet$ there exists a constant $C> 0$ such that we have $\mu_\rho(x)\leq C $ for any $x\in X$ and $\rho> 0$.

\noindent $\bullet$ we have $\lim_{\rho \to 0}\mu_\rho (x)= \nu \big(\Theta_{h_L}(L), x\big)$ for each $x\in X$.

\medskip
\noindent The discussion in the paragraph D still applies, but in the actual
circumstances we only obtain
$$\dbar \big(\sqrt {\eta_\varepsilon+ \lambda_\varepsilon}v_{\varepsilon, \rho}\big)+ \sqrt {\mu_\rho\eta_\varepsilon}w_{\varepsilon, \rho}^\prime+
\dbar^\star(\sqrt \eta_\varepsilon w_{\varepsilon, \rho}^{\prime\prime})= 
\dbar \Big(\theta\Big({|s|^2\over \varepsilon^2}\Big)U_\infty\Big)\leqno(38)$$
together with estimates for $v_{\varepsilon, \rho}, w_{\varepsilon, \rho}^{\prime}$
and $w_{\varepsilon, \rho}^{\prime\prime}$ similar to (26). The adjoint
operator $\dbar ^\star$ is with respect to the metric $\varphi_{L, \rho}+ \log |s|^2$ on $X\setminus Z$. 

Since we want to use the fact that the metric $h_{L, \rho}$ is non-singular,
we have to take first $\varepsilon\to 0$, and then $\rho\to 0$ in the relation 
above. The trouble comes from the presence of the factor 
$\dbar^\star(\sqrt \eta_\varepsilon w_{\varepsilon, \rho}^{\prime\prime}):$
in order to remove it, one can project on the kernel of the $\dbar$
(similar ideas were used in [2]). Indeed, the arguments invoked in [2] would give the desired result,
provided that 
one can obtain estimates for the $L^2$ norm of 
$$\sqrt {\mu_\rho\eta_\varepsilon}w_{\varepsilon, \rho}^\prime$$
and not just for the $L^2$ norm of 
$w_{\varepsilon, \rho}^\prime$ as it is the case here. The puzzling fact is that 
if one tries to obtain these estimates by incorporating the factor
$$\log\big(\log{1\over \varepsilon^2+ |s|^2}\big)$$
in the expression of the metric of $L$, then the curvature hypothesis needed are basically unchanged (which is very good!), but the right hand side term
of (35) becomes
$$\int_X\theta^\prime \Big({{\vert s\vert^2}\over {\varepsilon ^2}}\Big)^2\log {1\over \varepsilon^2+ |s|^2}{{|U_\infty|^2}\over {|s|^2}}
e^{-\varphi_{L, \rho}-\varphi_E}$$
which is unbounded as $\varepsilon\to 0$. In conclusion, we are one estimate far from the general case.


\section{\S G An open problem}

\vskip 7pt

\noindent The next  problem
was formulated by Y.-T. Siu in his seminal article 
[17].

\claim Conjecture {\rm (Siu)}|
 Let $\pi:{\cal X}\to \Delta$ be a
K\"ahler family over the unit disk.
Then any section of $\displaystyle mK_{{\cal X}_0}$ extends to
${\cal X}$.
\endclaim

\medskip It is rather safe to predict that the Ohsawa-Takegoshi theorem
in the K\"ahler setting will be necessary for the solution of this 
problem: this is a very strong motivation for our research in this article.
Given the proof of the conjecture above in the case $X$ projective, the hypothesis 
$$u\in H^{0}\big(Z,(K_{Z}+L_{|Z})\otimes 
{\cal I}_+(h_{L}|_{Z})\big)$$
is not as restrictive as it looks: the construction of the metric for the bundle 
$L$ in [17] shows that the initial section is in fact {\sl bounded}
with respect to this metric. A slight drawback in the statement 0.1 is
the hypothesis $\cC$, concerning the behavior of the singularities of 
$h_L$; we refer to [16] for a discussion of this property, in the context of
metrics with logarithmic singularities.

\section{References}

\bigskip

{\eightpoint

\bibitem [1]&Berndtsson, B.:& On the Ohsawa-Takegoshi extension theorem;& Ann.\ Inst.\ Fourier (1996)&

\bibitem [2]&Demailly, J.-P.:&Estimations $L^2$ pour l'op\'erateur d-bar d'un fibr\'e vectoriel holomorphe semi-positif au-dessus d'une vari\'et\'e k\"ahl\'erienne compl\`ete, &\ Ann. Sci. Ecole Norm. Sup. 4e S\'er. 15 (1982) 457-511&

\bibitem [3]&Demailly, J.-P.:&Regularization of closed positive currents of type (1,1) by the flow of a Chern connection,&\ Actes du Colloque en l'honneur de 
P.~Dolbeault (Juin 1992), \'edit\'e par H. Skoda et J.M. Tr\'epreau, Aspects of Mathematics, Vol. E 26 (1994) 105-126&

\bibitem [4]&Demailly, J.-P.:&  On the Ohsawa-Takegoshi-Manivel  
extension theorem;& Proceedings of the Conference in honour of the 85th birthday of Pierre Lelong, 
Paris, September 1997, Progress in Mathematics, Birkauser, 1999&

\bibitem [5]&Demailly, J.-P, Peternel, Th.:&
A Kawamata-Viehweg Vanishing Theorem on Compact K\"ahler Manifolds
&\ J. Differential Geom. Volume 63, Number 2 (2003), 231-277&

\bibitem [6]&H\"ormander, L.:&
$L^2$ estimates and existence theorems for the $\dbar$ operator;&\ Acta Mathematica
Volume 113, Number 1, 89-152&

\bibitem [7]&Kim, D.:&$L^2$ extension of adjoint line bundle sections,;&
arXiv:0802.3189, to appear in Ann.\ Inst.\ Fourier&

\bibitem [8]&Koziarz, V.:&\ Extensions with estimates of cohomology classes;&  preprint 
available on the web page of the author, 2009&

 \bibitem [9]&Manivel, L.:& Un th\'eor\`eme de prolongement L2 de sections holomorphes d'un fibr\'e hermitien;& Math.\ Zeitschrift {\bf 212} (1993), 107-122&
 
 \bibitem [10]&McNeal, J., Varolin, D.:& Analytic inversion of adjunction;&\ 
 Ann.\ Inst.\ Fourier (Grenoble)  {\bf 57}  (2007),  no. 3, 703--718&
 
\bibitem [11]&Morrow, J., Kodaira, K.:&\ Complex manifolds&\ AMS Chelsea
  Publications&
 
\bibitem [12]&Ohsawa, T., Takegoshi, K.\ :& On the extension of $L^2$
holomorphic functions;& Math.\ Z.,
{\bf 195} (1987), 197--204&

 \bibitem [13]&Ohsawa, T\ :& A precise $L\sp 2$ division theorem;&
 Complex geometry (G\"ottingen, 2000), 185--191, Springer, Berlin, 2002. &
 
 \bibitem [14]&Ohsawa, T.\ :& On the extension of $L\sp 2$ holomorphic functions. VI. A limiting case;&
Contemp.\ Math. \ ,  {\bf 332} (2003)
Amer. Math. Soc., Providence, RI&

 \bibitem [15]&Pali, N.:&On the extension of twisted holomorphic sections of singular hermitian line bundles&\ arXiv:0807.3535&
 
 \bibitem [16]&Phong, D. H., Sturm, J.&\ Algebraic estimates, stability of
local zeta functions, and uniform
estimates for distribution functions&\ Annals of Mathematics, 152 (2000)&
 
\bibitem [17]&Siu, Y.-T.:&\ Extension of twisted pluricanonical sections with plurisubharmonic weight and invariance of semipositively twisted plurigenera for manifolds not necessarily of general type;& Complex geometry (G\"ottingen, {\bf 2000}),  223--277, Springer, Berlin, 2002&

\bibitem [18]&Siu, Y.-T.:&\ Multiplier ideal sheaves in complex and algebraic geometry;& Sci.\ China Ser.  {\bf A 48}, 2005&

\bibitem [19]&Varolin, D.:&\ A Takayama-type extension theorem;&  math.CV/0607323, Comp.\ Math&

}
\medskip

\it IECN, Universit\'e Henri Poincar\'e Nancy 1 B.P. 70239, 54506 Vandoeuvre-l\`es-Nancy Cedex, France\\
\it E-mail adress: yi@iecn.u-nancy.fr

\end